\documentclass[letterpaper, 10 pt, conference]{ieeeconf}  % Comment this line out if you need a4paper

\IEEEoverridecommandlockouts                              % This command is only needed if 
\usepackage[utf8]{inputenc}
\usepackage{amssymb}
\usepackage{amsmath}
\usepackage{bm}
\usepackage{hyperref}
\usepackage{tikz}
\usepackage{tkz-euclide}
\usepackage{chngcntr}
\usepackage{verbatim}
\usepackage{mathtools}
\usepackage{esvect}
\usepackage{subfiles}
\usepackage{subfig}
\usepackage{cite}
\usepackage{color}
\usepackage{xspace}
\usepackage[noend]{algpseudocode}
\RequirePackage{algorithm}

\newtheorem{theorem}{Theorem}
\newtheorem{lemma}[theorem]{Lemma}

\newtheorem{assumption}{Assumption}
\newtheorem{definition}{Definition}
\newcommand{\norm}[1]{\left\|#1\right\|}
\newcommand{\Oc}{\mathcal{O}}
\newcommand{\reals}{\mathbb{R}}
\newcommand{\inner}[1]{\langle #1 \rangle}
\newcommand{\grad}{\nabla}
\newcommand{\vectwo}[2]{\begin{bmatrix} #1 \\ #2 \end{bmatrix}}

\newcommand{\cE}{\mathcal{E}}
\newcommand{\cO}{\mathcal{O}}

\title{\LARGE \bf
 Acceleration in First Order Quasi-strongly Convex \\
 Optimization by ODE Discretization
}

\author{Jingzhao Zhang, Suvrit Sra and Ali Jadbabaie % <-this % stops a space
	\thanks{The authors are with the Laboratory for Information and Decision Systems (LIDS), and the Institute for Data, Systems, and Society (IDSS),
		Massachusetts Institute of Technology, 77 Massachusetts Ave, Cambridge, MA 02139
		{\tt \{jzhzhang, suvrit, jadbabai\}@mit.edu}}%
	\thanks{*This research was supported in part by DARPA Lagrange.}
}

\begin{document}

	\maketitle
	%%%%%%%%%%%%%%%%%%%%%%%%%%%%%%%%%%%%%%%%%%%%%%%%%%%%%%%%%%%%%%%%%%%%%%%%%%%%%%%%
	\begin{abstract}
		
		\textcolor{black}{We study gradient-based optimization methods obtained by direct Runge-Kutta discretization of the ordinary differential equation (ODE) describing the movement of a heavy-ball under constant friction coefficient. When the function is high order smooth and strongly convex, we show that directly simulating the ODE with known numerical integrators achieve acceleration in a nontrivial neighborhood of the optimal solution. In particular, the neighborhood can grow larger as the condition number of the function increases. Furthermore, our results also hold for nonconvex but quasi-strongly convex objectives. We provide numerical experiments that verify the theoretical rates predicted by our results.}
		
	\end{abstract}

	%%% TeX-master: "../root"
\section{Introduction}
In this paper, we study accelerated first-order optimization algorithms for the problem
\begin{equation}
\label{main_problem}
\min_{x \in \reals^d}\quad f(x),
\end{equation}
where $f$ satisfies the following properties.
\begin{assumption} \label{assump:quasi-sc}
	Denote $x^*$ the unique optimal solution of $f$. Assume function $f$ is quasi-strongly convex, i.e. $\forall x$,
	\begin{align*} 
		f(x^*) \ge f(x) + \inner{\nabla f(x), x^* - x} + \tfrac{\mu}{2}\|x - x^*\|^2.
	\end{align*}
\end{assumption}
\begin{assumption}\label{assump:smooth}
	Function $f$ is $L-$smooth, i.e. $\forall x, y,$
	\begin{align*}
		\|\nabla f(x) - \nabla f(y)\| \le L\|x - y\|.
	\end{align*}
\end{assumption}

Assumption~\ref{assump:quasi-sc} is implied by the standard strong-convexity assumption and is weaker as the inequality only needs to hold when one of the two points is the optimal solution. Further, Assumption~\ref{assump:quasi-sc} does not imply convexity (see \cite{necoara2018linear}). A classical method for solving~\eqref{main_problem} is gradient descent (GD), which requires $T = \Oc(\frac{L}{\mu}\log(\tfrac{f(x_0) - f(x^*)}{\epsilon})$ iterations to achieve $\epsilon$ accuracy ($f(x_T) - f(x^*) \le \epsilon$). With the additional global convexity assumption, it is shown in \cite{necoara2018linear} that Nesterov's accelerated gradient (NAG) method~\cite{nesterov-smooth-acceleration} matches the oracle lower bound of $O(\sqrt{\frac{L}{\mu}}\log(\tfrac{f(x_0) - f(x^*)}{\epsilon})$~\cite{nemirovskii-lowerbound}. However, up to our knowledge, no acceleration in terms of condition number dependency is achieved without the convexity assumption.

In this work, we show that with the higher order smoothness assumption, known Runge-Kutta discretization of the standard heavy-ball ODE~\eqref{eq:hb-ode} achieves better dependency of $\frac{L}{\mu}$ compared to gradient descent in a neighborhood of the optimal points. Unlike local results, we show that the size of the neighborhood can potentially increases with the condition number $\frac{L}{\mu}$. For more details, please refer to Theorem~\ref{thm:main}.
\vspace{-0.2cm}
\section{Related work}
\vspace{-0.2cm}
Acceleration of first order optimization algorithms has long been studied. Polyak proposed the heavy ball method~\cite{polyak1964some} that is provably faster than gradient descent for optimizing quadratic objectives.  Nesterov later designed acceleration methods for general smooth convex objectives~\cite{nesterov-smooth-acceleration}. It was recently proven in~\cite{necoara2018linear} that this algorithm also achieves optimal rate under global convex and quasi-strongly convex conditions. However, ever since its introduction, acceleration has remained somewhat mysterious, especially because Nesterov's original derivation relies on elegant but unintuitive algebraic arguments. This lack of understanding has spurred a variety of recent attempts to uncover the rationale behind the phenomenon of acceleration~\cite{allen-linear, bubeck-geometric, lessard-iqc, hu-dissipativity, scieur-regularized, fazlyab-dynamical}.

Our approach follows the sequence of work on the ODE  interpretation of optimization algorithms. The ODE interpretation has long been studied, for example in~\cite{alvarez2000minimizing, attouch2000heavy, bruck1975asymptotic, attouch1996dynamical}. These work analyzed the asymptotic behavior of dissipative dynamical systems. The more recent sequence of work on the ordinary ODE interpretation of the acceleration took more interests in the non-asymptotic analysis and discretized algorithms. This line of research starts with~\cite{su-differential}, who showed that the continuous limit of NAG is a second order ODE describing a physical system with vanishing friction.  Later, \cite{wibisono-variational, xu2018accelerated, krichene-accelerated,  francca2018admm, barakat2018convergence,wilson2019accelerating, maddison2018hamiltonian} generalized this idea to study the continuous limit of other first order algorithms such as mirror descent, ADMM and rescaled gradient descent. \cite{shi2018acceleration, zhang2018achieving, zhang2018direct, betancourt2018symplectic, shi2018understanding} studied the discretization of ODE and aim to design accelerated algorithms with non-asymptotic convergence guarantee. 

Most of the literature requires knowing Nesterov's method beforehand to study acceleration. Our contribution is to show that acceleration can be achieved starting from the heavy-ball ODE without knowing NAG. The most relevant work to ours would be \cite{zhang2018direct,shi2018acceleration}. We highlight two differences.  First, unlike \cite{shi2018acceleration} which discretizes a high order ODE approximation of NAG, we discretize an ODE that describes the heavy ball mechanical system. Compared to \cite{shi2018acceleration}, our approach can produce accelerated algorithms without knowing NAG first. Second, we focus on quasi-strongly convex problems, which is different from the setting of both work.

	%%% TeX-master: "../root"
\section{Preliminaries}
Before presenting the main result, we will go over some tools used in our algorithms and analysis. The first one is a class of numerical discretization algorithms named Runge-Kutta integrators. The second one is a notation system that simplifies the expression for high order derivatives of autonomous dynamical systems.

\subsection{Runge-Kutta Integration}
We briefly recall \emph{explicit} Runge-Kutta (RK) integrators used in our work. For a more in depth discussion please see the textbook~\cite{hairer-textbook}. 
\begin{definition} \label{def:rk}
	Given a dynamical system $\dot{y} = F(y)$, let the current point be $y_0$ and the step size be $h$. An explicit $S$ stage Runge-Kutta method generates the next step via the following update:
	\begin{align*}\label{eq:runge-kutta}    
	g_i &= y_0 + h \sum_{j=1}^{i-1} a_{ij} F(g_j), \\
	\Phi_h(y_0) &= y_0 + h \sum_{i=1}^{S} b_i F(g_i),
	\end{align*}
	where $a_{ij}$ and $b_i$ are suitable coefficients defined by the integrator;  $\Phi_h(y_0)$ is the estimation of the state after time step $h$, while $g_i$ (for $i=1,\ldots,S$) are a few neighboring points where the gradient information $F(g_i)$ is evaluated. %The summation limit $S$ in~\eqref{eq:runge-kutta} is named the number of stages.
\end{definition}
By combining the gradients at several evaluation points, the integrator can achieve higher precision by matching up Taylor expansion coefficients. Let $\varphi_h(y_0)$ be the true solution to the ODE with initial condition $y_0$. Then we define the order of an integrator below. 
\begin{definition} \label{def:order}
we say that an integrator $\Phi_h(y_0)$ has order $s$ if its \emph{discretization error} shrinks as
\begin{equation}
\label{eq:disc-error}
\|\Phi_h(y_0) - \varphi_h(y_0)\| = O(h^{s+1}), \qquad\text{as}\ h\to 0.
\end{equation}
\end{definition}
In general, RK methods offer a powerful class of numerical integrators, encompassing several basic schemes. The \emph{explicit Euler's} method defined by $\Phi_h(y_0) = y_0 + hF(y_0)$ is an explicit RK method of order 1, while the \emph{midpoint} method $\Phi_h(y_0) = y_0 + hF(y_0+\tfrac{h}{2}F(y_0))$ is of order 2. Some high-order RK methods are summarized in \cite{highorderRK}. An order 4 RK method requires 4 stages, i.e., 4 gradient evaluations, while an order 9 method requires 16 stages.

\subsection{Elementary differentials}  \label{sec:elementary}
\begin{figure}[t]
	\centering
	\includegraphics[width=0.4\textwidth]{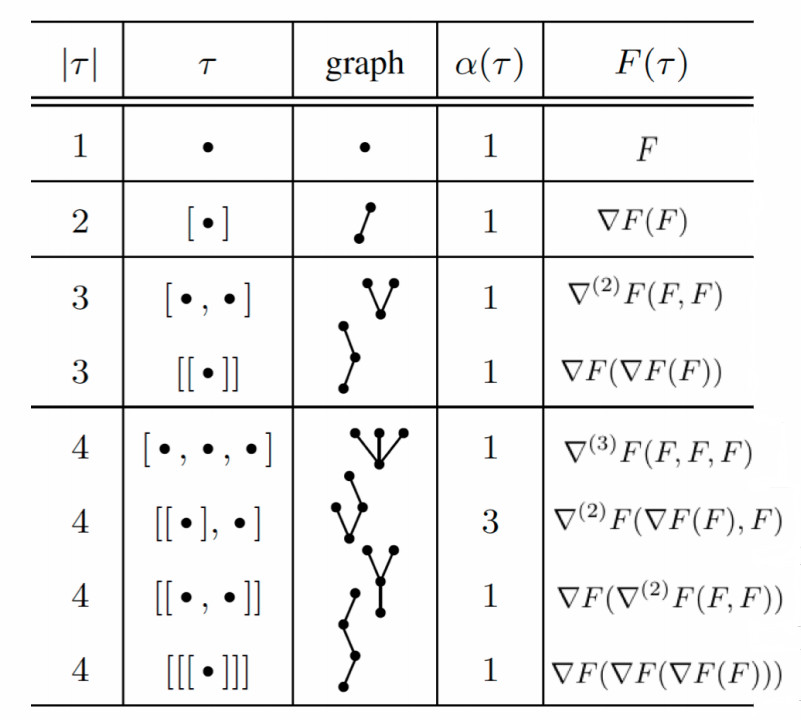}
	\caption{An illustration from \cite{hairer-textbook}. Example tree structures and corresponding function derivatives.}
	\label{fig:bseries}\vspace{-0.5cm}
\end{figure}

We briefly summarize some key results on elementary differentials from \cite{hairer-textbook}. For more details, please refer to chapter 3 of the book. Given a dynamical system 
$$\dot{y} = F(y),$$
we want to find a convenient way to express and compute its higher order derivatives. To do this, let $\tau$ denote a tree structure as illustrated in Figure~\ref{fig:bseries}. $|\tau|$ is the number of nodes in $\tau$. Then we can adopt the following notations as in \cite{hairer-textbook}.

\begin{definition} \label{def:elementary}
	For a tree $\tau$, the elementary differential is a mapping $F(\tau): \reals^d \to \reals^d$, defined recursively by $F(\bullet)(y) = F(y)$ and 
	$$F(\tau)(y) = \nabla^{(m)}F(y)[F(\tau_1)(y), ..., F(\tau_m)(y)],$$
	for $\tau = [\tau_1, ..., \tau_m]$. Notice that $\sum_{i=1}^m |\tau_i| = |\tau| - 1$.
\end{definition}

Some examples are shown in Figure~\ref{fig:bseries}. With this notation, the following results from \cite{hairer-textbook} Chapter 3.1 hold. The proof follows by recursively applying the product rule.
\begin{lemma} \label{lemma:high-derivative-solution}
	The qth order derivative of the exact solution to $\dot{y} = F(y)$ is given by
	$$\frac{d^q y(t)}{dt^q}|_{t=t_c} = y^{(q)}(t_c) = F^{(q-1)}(y_c) = \sum_{|\tau|=q}\alpha(\tau)F(\tau)(y_c),$$
	for $y(t_c) = y_c$. $\alpha(\tau)$ is a positive integer determined by $\tau$ and counts the number of occurrences of the tree pattern $\tau$.
\end{lemma}
The next result is obtained by general Leibniz rule. The expression for $\frac{\partial^q F(g_i)}{\partial h^q}$ can be calculated the same way as in Lemma \ref{lemma:high-derivative-solution}.
\begin{lemma} \label{lemma:high-derivative-numerical}
	For a Runge-Kutta method defined in definition ~\ref{def:rk}, if $F$ is $q_{th}$ differentiable, then
	\begin{align}\label{eq:Leibniz}
	\frac{d^q \Phi_h(y_c)}{d h^q} = \sum_{i \le S} b_i (h \frac{d^q F(g_i)}{d h^q} + q\frac{d^{q-1} F(g_i)}{d h^q}),
	\end{align}
	where $\frac{\partial^q F(g_i)}{\partial h^q}$ has the same structure as $F^{(q)}(y)$ in lemma~\ref{lemma:high-derivative-solution}, except that we need to replace all $F$ in the expression by 
	$\frac{\partial g_i}{\partial h}$ and all $\nabla^{(n)}F(y)$ by $\nabla^{(n)}F(g_i)$.
\end{lemma}

	%%% TeX-master: "../root"

\section{Main results}
In this section, we will propose a family of algorithms generated via ODE discretization. We will also provide theoretical convergence guarantee and go over the proofs. 

\subsection{Algorithm}
In this subsection, we introduce a second-order dynamical system and use explicit RK integrators to generate iterates that converge to the optimal solution at a rate faster than gradient descent ($\Oc(\frac{L}{\mu}\log(\tfrac{f(x_0) - f(x^*)}{\epsilon})$). We start with the second-oder ODE,
\begin{align} \label{eq:hb-ode}
\ddot{x}(t)   + 2\dot{x}(t)  + \frac{\nabla f(x)}{\mu}  = 0.
\end{align}
This ODE has very nice physical interpretation. It describes the mechanics of a heavy ball moving on the surface defined by $f$. Intuitively, due to friction, the heavy ball will stop at the unique optimal of $f$. This second order ODE can be written as an autonomous dynamical system.
\begin{align}\label{eq:dynamic}
\dot{y} = \vectwo{\dot{v}}{\dot{x}} = F(y) = \begin{bmatrix} -2v -\grad f(x)/\mu\\ v \end{bmatrix}.
\end{align}
As the stability of discretization is determined by the largest numerical error along any coordinate, we would like to \textbf{ balance the Lipschitz constant of $x$ and $v$}. This step is crucial to achieving acceleration. Particularly, we make a change of variable $w = v/\sqrt{Q}$, where $Q = \frac{L}{\mu}$ is the condition number. Hence with some abuse of notations, we can rewrite the previous system in the following form,
\begin{align} \label{eq:re-dynamic}
\dot{y} = \vectwo{\dot{w}}{\dot{x}} = F(y) = \begin{bmatrix} -2w -\frac{\grad f(x)}{\mu\sqrt{Q}}\\ \sqrt{Q}w \end{bmatrix}.
\end{align}
The algorithm we propose is simply to discretize the dynamical system in~\eqref{eq:re-dynamic} with an order-$s$ explicit Runge-Kutta integrator. It is summarized in Algorithm 1.

\begin{algorithm}[t] \label{algorithm_main}
	Algorithm 1: Input($f, x_0, L, \mu, s, N$)\\
	\Comment{Constants $\mu, L$ are the same as in Assumptions}
	\begin{algorithmic}[1]
		\State Set the initial state $y_0 = [\vec{0}; x_0]\in \reals^{2d}$
		\State Set step size h = $C/N^{\frac{1}{s+1}}$. \Comment{C is determined by $\mu, L, s, x_0$}
		\State $x_N \gets \text{Order-s-Runge-Kutta-Integrator}(F, y_0, N, h)$ \Comment{F is defined in equation~\eqref{eq:re-dynamic}}
		\State \textbf{return} $x_N$
	\end{algorithmic}
\end{algorithm}

\subsection{Convergence analysis}
In order to utilize the order conditions of Runge-Kutta integrators, we make the following assumptions on high order smoothness. Denote $\| M \| $ the operator norm of an $n_{th}$ order tensor 
$M: \reals^{\overbrace{d\times...\times d}^n} \to \reals.$
In other word,
\begin{align*}
\|M\| = \max \{ M[u_1, ..., u_n] | u_i \in \reals^d, \|u_i\|_2 = 1 \}
\end{align*}
Then we assume,
\begin{assumption}\label{assump:high-smooth}
	The high order derivatives of $f$ exist and are bounded up to order $s+1$. In other words,
	$\|\grad^{(n)} f(x) \| \le L$, $n = 2, ..., s+1$
\end{assumption}

For this nonlinear dynamical system, we can define the following Lyapunov function to help us prove convergence in both continuous time and discrete time.
\begin{align}\label{eq:e0}
\cE(y) = &\cE([w; x]) = 2(f(x) - f(x^*))/\mu + \frac{Q}{2}\|w\|^2 + \nonumber \\
&\frac{1}{2}\|x + \sqrt{Q}w - x^*\|^2.
\end{align}
The stability of the continuous system is justified by the following Lemma.
\begin{lemma} \label{lemma:cts-exp-conv}
	For the dynamical system defined in~\eqref{eq:re-dynamic}, the Lyapunov function converges exponentially along the trajectory with rate
	$$\dot{\cE}(y) \le -\frac{1}{2}\cE(y).$$
\end{lemma}
\begin{proof}[of Lemma ~\ref{lemma:cts-exp-conv}]
\begin{align*}
\dot{\cE}(y) &= \tfrac{2}{\mu}\inner{\grad f(x), \sqrt{Q}w} + Q \inner{w, -2w-\grad f(x)/(\sqrt{Q}\mu)} \\
&+  \inner{x + \sqrt{Q}w - x^*, - \sqrt{Q}w - \tfrac{\grad f (x)}{\mu}} \\
&= -2Q\|w\|^2 - \inner{x - x^*, \grad f(x)}/\mu \\
& - \inner{x + \sqrt{Q}w - x^*, \sqrt{Q}w}\\
&\le -2Q\|w\|^2 - (f(x) - f(x^*) + \mu/2\|x - x^*\|^2)/\mu \\
& - \tfrac{1}{2} (\|x + \sqrt{Q}w - x^*\|^2 + \|\sqrt{Q}w\|^2 - \|x- x^*\|^2) \\
& \le -\cE(y)/2.
\end{align*}
The first inequality follows from Assumption~\ref{assump:quasi-sc} and the fact that $ \inner{a, b} = \frac{1}{2}(\|a\|^2 + \|b\|^2 - \|a-b\|^2)$.
\end{proof}

Then we provide the following convergence guarantee.
\begin{theorem} \label{thm:main}
	Under the Assumptions~\ref{assump:quasi-sc}~\ref{assump:high-smooth}, we discretize the dynamical system~\eqref{eq:re-dynamic} with an order-$s$ Runge-Kutta integrator. By setting $h \le \min\{\cE(y_0)^{-\gamma}, Q^{-\gamma} \}\frac{1}{2c(s)^{\frac{1}{s+1}}}, \gamma = \frac{s+3}{2(s+1)}$ for some constant $c(s)$ determined by the order $s$, we have,
	$$\cE(y_N) \le (1 - h/4)^N\cE(y_0). $$
	In other word, assuming the integrator has stage $S$, the number of gradients evaluated to get $\epsilon$ suboptimality is $\cO(S\max\{\cE(y_0), Q\}^\gamma \log(\frac{\cE(y_0)}{\mu\epsilon})).$ 
	
\end{theorem}

We would like to make two remarks about this result. First, for general quasi-strongly convex functions (with the exception of quadratics), the smoothness constant $L$ would grow with $\mathcal{E}_0$. Second, when the initial point is selected such that $\mathcal{E}_0 \le \frac{L}{\mu}$, the algorithm achieves acceleration and approaches the theoretical lower bound as the order of the integrator $s$ increases.

\subsection{Proof of main theorem}
The high level idea is simple. Based on the convergence result of Lemma~\ref{lemma:cts-exp-conv}, we get exponential convergence as long as the discretization error stays bounded. Then the convergence rate would depend on the largest stable discretization step. We postpone the technical details to the next subsection and abstract out the main proof below.

Recall that $\varphi_h(y_0)$ is the true solution to the ODE with initial condition $y_0$, while $\Phi_h(y_0)$ is the numerical solution generated by the integrator.  We rewrite 
\begin{align} \label{eq:def-yh}
	y(h) = \varphi_h(y), \quad g_S(h) = \Phi_h(y),
\end{align}
to emphasize the explicit dependency on $h$.
By the order condition of the Runge-Kutta integrator, we know that the first $s$ terms in the Taylor expansions of $y_0(h)$ and $g_S(h)$ evaluated at $h=0$ are the same. By Taylor's theorem and triangle inequality, we know that
\begin{align*}
|\cE(g_S(h))& - \cE(y(h))| \le \\
 &h^{s+1} \max_{0 \le \delta \le h} (|\frac{d^{s+1}}{dh^{s+1}} \cE(g_S(\delta)) | + |\frac{d^{s+1}}{dh^{s+1}} \cE(y(\delta)) | ).
\end{align*}
Since $\cE(y(h)) \le (1 - h/2)\cE(y)$, we have
\begin{align*}
\cE(g_S(h) \le &(1 - h/2)\cE(y) \\
+ & h^{s+1} \max_{0 \le \delta \le h} (|\frac{d^{s+1}}{dh^{s+1}} \cE(g_S(\delta)) | + |\frac{d^{s+1}}{dh^{s+1}} \cE(y(\delta)) | ) \\
\le &(1 - h/2)\cE(y) \\
+ & c(s)h^{s+1} \sum_{i=1}^{s}Q^{(i+1)/2}\cE(y)^{(s-i+2)/2}\\
\le &(1 - h/4)\cE(y_0).
\end{align*}
The second inequality follows by Lemma~\ref{lemma:high-order-dE-dt} and Lemma~\ref{lemma:high-order-dE-dh} which bound the two derivative terms respectively. The last inequality follows by the choice of $h$.

\subsection{Technical lemmas}
In this section, we bound the high order derivatives of the dynamical system and the discretization algorithm as a function of time. Before presenting a few key lemmas, we  define some quantities and frequently used inequalities. 

\begin{lemma} The norm of the dynamics is bounded by the Lyapunov function, 
$$\|F(y)\| \le 5\sqrt{\cE(y)}.$$
\end{lemma}
\begin{proof}
By~\eqref{eq:re-dynamic},
	\begin{align*}
		\|F(y)\|^2 &=\left \| \begin{bmatrix} -2w -\frac{\grad f(x)}{\mu\sqrt{Q}}\\ \sqrt{Q}w \end{bmatrix} \right\|^2 \\
		&\le 2\|\frac{\grad f(x)}{\mu\sqrt{Q}}\|^2 + 8\|w\|^2 + Q\|w\|^2.
	\end{align*}
	By Lipschitz  continuity of the gradient and optimality of $f(x^*)$, we get that $\|\grad f(x)\|^2 \le L(f(x) - f(x^*))$ (see \cite{nesterov2013introductory}). Therefore, we have
	\begin{align*}
	\|\frac{\grad f(x)}{\mu\sqrt{Q}}\|^2 &= \|\grad f(x)\|^2/(\mu L)   \le (f(x) - f(x^*))/\mu.
	\end{align*}
	Substitute this in the inequality above and we proved the lemma.
\end{proof}

Next, we bound the high order derivatives of the system dynamics $F$ defined in \eqref{eq:re-dynamic}.  
\begin{lemma} \label{lemma:grad-F}
	$\|\grad^{(n)} F(y)\| \le 4\sqrt{Q}, \forall n \ge 1$.
\end{lemma}
 \begin{proof}
 	This follows by Assumption \ref{assump:high-smooth} and computing $\|\grad^{(n)} F(y)\|$ explicitly.
 \end{proof}
 
 Next, notice
\begin{align*}
\grad \cE(y) = \begin{bmatrix} 
Qw + \sqrt{Q}(x + \sqrt{Q}w - x^*)\\
2\grad f(x)/\mu + (x + \sqrt{Q}w - x^*) \end{bmatrix}.
\end{align*}
From this expression, we can tell that $\|\grad \cE(y)\|^2 \le 4Q\cE(y)$. Further, it's easy to verify that $\|\grad^{(n)} \cE(y)\| \le 5Q, \forall n > 2$. 

The above inequalities serve as the base case of the Lemmas proven below with induction and will be used repeatedly.

In the rest of this section, we let $F$ denote the function defined in \eqref{eq:re-dynamic} and $F(\tau )$ denotes the differentials defined in section~\ref{sec:elementary}. Then the next Lemma follows. For conciseness, we write $\cE$ as a shorthand for $\cE(y)$.

\begin{lemma}\label{lemma:F-tau}
	The vector norms of elementary differentials $F(\tau)$ evaluated at $y$ can be bounded by the value of the Lyapunov function $\mathcal{E}(y)$. More precisely, $\forall k \ge 2, \sum_{|\tau| = k}  \|F(\tau)\| \le  c(k) \sum_{i = 1}^{k-1} Q^{\tfrac{k-i}{2}}\cE^{\tfrac{i}{2}} $.
\end{lemma}

\begin{proof}
	We already know $\|F(\bullet)\| = \|F(y)\| \le  5\sqrt{\cE} $ and $\norm{\grad F(y)} \le 4\sqrt{Q}$. Hence the case when $k=2$ follows by Cauchy-Schwartz. We prove the Lemma by induction and assume that for any tree $|\tau| = k \le l $,
	\begin{align*}
	\sum_{|\tau| = k} \| F(\tau) \| \le  c(k) \sum_{i = 1}^{k-1} Q^{\tfrac{k-i}{2}}\cE^{\tfrac{i}{2}}.
	\end{align*}
	By definition of tree structures, for $\tau' = [\tau_1, ..., \tau_m]$ with $|\tau'| = l+1, m \ge 1$, 
	We have 
	\begin{align*}
	\|F(\tau')\| \le \| \grad^{(m)} F \| \prod_{i=1}^m \|F(\tau_i)\|.
	\end{align*}
	Therefore, by the inductive assumption for some absolute onstant $c(l)$ determined by $l$,
	\begin{align*}
	\sum_{|\tau'| = l+1} \|F(\tau')\| \le c(l+1) \sum_{i = 1}^{l} Q^{\tfrac{l+1-i}{2}}\cE^{\tfrac{i}{2}}.
	\end{align*}
	
\end{proof}

With the above Lemma, we are able to show that the high order derivatives of the dynamical system are bounded. 
\begin{lemma}\label{lemma:high-order-dy-dt}
Recall the definition of $y(h)$ in \eqref{eq:def-yh}.
$$\forall k \ge 2, \norm{\frac{d^k}{dh^k} y(h) }  \le  c(k) \sum_{i = 1}^{k-1} Q^{\tfrac{k-i}{2}}\cE^{\tfrac{i}{2}}.$$
\end{lemma}
\begin{proof}
	This follows directly by Lemma \ref{lemma:high-derivative-solution} and \ref{lemma:F-tau}.
\end{proof}

We have already bounded the high order derivatives of the continuous trajectory. However, notice that the bounds are point-wise in time. We would like to bound the variation of the Lyapunov function for a short time span. 
\begin{lemma} \label{lemma:lyap-small-h}
	If we set $h \le \frac{1}{10\sqrt{Q}}$,  then $\cE(y(h)) \le 3\cE(y(0)).$
\end{lemma}
\begin{proof}
	By Taylor expansion, we get
	\begin{align*}
	\cE(y(h)) &= \cE(y) + h\frac{d}{dh}\cE(y(h))|_{h=0}\\
	&+ \frac{h^2}{2}\int_0^h \frac{d^2}{dh^2} \cE(y(h)) d\tau.
	\end{align*}
	Since $ \dot{\cE}(y(0)) \le -\frac{1}{2}\cE(y) $, we have
	\begin{align*}
	\cE(y(h))& \le (1-h/2)\cE(y) + \frac{h^2}{2}\int_0^h \frac{d^2}{d\tau^2} \cE(y(\tau)) d\tau.
	\end{align*}
	
	We know that 
	\begin{align*}
	\frac{d^2}{d\tau^2} \cE(y(\tau))) = \grad^{(2)} \cE(y(\tau))[F(y(\tau)), F(y(\tau))] \\
	+ \grad \cE(y(\tau))^T \grad F(y(\tau)) F(y(\tau)).
	\end{align*}
	Notice also that from previous Lemmas,
	\begin{align*}
	&\|\grad^{(2)} \cE(y(\tau))\| \le 5Q, \quad
	\|F(y(\tau))\| \le 5\sqrt{\cE(y(\tau))}, \\
	&\| \grad F(y(\tau))) \| \le 4\sqrt{Q}, \quad
	\|\grad \cE(y(\tau))\| \le 2\sqrt{Q\cE(y(\tau)) }.
	\end{align*}
	This implies 
	\begin{align*}
		\cE(y(h))& \le (1-h/2)\cE(y) + 100h^2Q \int_0^h \cE(y(\tau)) d\tau \\
		&\le \cE(y) +  \int_0^h \cE(y(\tau)) \le e^h \cE(y) \le 3\cE(y).
	\end{align*}
	The last inequality follows by $h\le 1$.
	
\end{proof}

We can then bound the derivatives of the Lyapunov function.

\begin{lemma}\label{lemma:high-order-dE-dt}
	$\forall k \ge 2, h \le 1/(2\sqrt{Q}),$
	$$\norm{\frac{d^k \cE(y(h))}{dh^k}} \le  c(k)\sum_{i=1}^{k-1}Q^{(i+1)/2}\cE^{(k-i+1)/2}(y(0)).$$
\end{lemma}
\begin{proof}
By chain rule, we have
	\begin{align*}
	\frac{d^n}{dt^n}\cE(y(t)) = \sum_{k_1, ..., k_n} \frac{n!}{k_1!k_2!...k_n!} \grad^{(k)} \cE(y) \prod_{i=1}^n (\frac{d^{i} y(t)}{dt^i}/i!)^{k_i}.
	\end{align*}
	where the sum is taken over $\{k_1, ..., k_n \in \mathbb{Z}_{\ge 0} | \sum_{i=1}^n ik_i = n\}$, and $k = \sum_{i=1}^n k_i$. Then the Lemma follows by Cauchy inequality, Lemma \ref{lemma:high-order-dy-dt}, Lemma \ref{lemma:lyap-small-h}, and $\|\grad^{(k)} \cE\| \le 5Q, \forall k \ge 2.$
	
\end{proof}

We have bounded the high order derivatives of the Lyapunov function along the true solution of the ODE. Next we will bound the counterpart for the points generated by discretization. Most of the steps would be almost the same. Hence we eliminate details due to limited space.

\begin{lemma} \label{lemma:dg-dh}
	Recall that $g_i(y)$ defined in Definition \ref{def:rk} is a function of the step size $h$. When $h \le \frac{1}{10\sqrt{Q}}$, $\forall i = 1, ...,  S$, we have $\|\frac{d}{dh}g_i(y)\| \le c\sqrt{\cE(y)}$ for some constant $c$ determined by the integration parameters $a_{ij}$ in Definition \ref{def:rk}. 
\end{lemma}
\begin{proof}
	Notice that $\frac{d}{dh}g_1(y) = a_{10}F(y)$. This proves the base case. The we prove by induction and the fact that
	\begin{align*}
	\frac{d}{dh}g_i(y) &= \sum_{j < i} a_{ij} F(g_j) + h\sum_{j < i} a_{ij} \grad F(g_j) \frac{d}{dh}g_j(y).
	\end{align*}
	The  claim follows by the inductive assumption and $\|\grad F(g_j)\| \le 4\sqrt{Q}$ (by Lemma ~\ref{lemma:grad-F}).
\end{proof}

With the above lemma, we can bound the Lyapunov of the neighboring points $g_i(y)$.
\begin{lemma}
	If we set $h \le \frac{1}{c\sqrt{Q}}$, where $c$ is determined by the numerical integrator, then $\cE(g_i(h)) \le 3\cE(y).$
\end{lemma}
\begin{proof}
The proof is almost the same as Lemma~\ref{lemma:lyap-small-h}.
\end{proof}
Then we bound the equivalent of Lemma~\ref{lemma:high-order-dE-dt}.

\begin{lemma} \label{lemma:high-order-dg-dh}
	$\forall i \le S, \exists c $ determined by $k$ and the integrator such that,  
	$$\|  \frac{d^k}{dh^k}g_i(y) \| \le  c(k, S) \sum_{i = 1}^{k-1} Q^{\tfrac{k-i}{2}}\cE^{\frac{i}{2}}.$$
\end{lemma}
\begin{proof}
The proof can be carried out in the same way as in Lemma~\ref{lemma:high-order-dy-dt}, except that we need to make use of Lemma~\ref{lemma:high-derivative-numerical} instead of Lemma~\ref{lemma:high-derivative-solution}. 
\end{proof}
Finally, we can bound the high order derivative of the Lyapunov function along the discretization solution.
\begin{lemma}\label{lemma:high-order-dE-dh}
	$\forall k \ge 2, h \le 1/(c\sqrt{Q}),$ where $c$ is determined by the integrator and $k$, 
	$$\norm{\frac{d^k \cE(g_S(h))}{dh^k}} \le  C\sum_{i=1}^{k-1}Q^{(i+1)/2}\cE^{(k-i+1)/2}(y(0)). $$
	The constant $C$ is determined by $k$ and the integrator. 
\end{lemma}
\begin{proof}
	The proof is exactly the same as in Lemma~\ref{lemma:high-order-dE-dt} except that instead of using Lemma \ref{lemma:high-order-dy-dt}, we use Lemma \ref{lemma:high-order-dg-dh} to bound the high order derivatives of the trajectories.
\end{proof}

	%%% TeX-master: "../root"
\section{Numerical experiments}
\begin{figure}[t]
	\centering
	\includegraphics[width=0.38\textwidth]{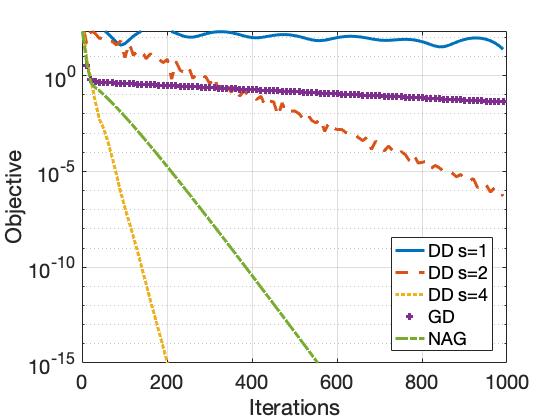}
	\caption{Convergence paths for GD, NAG and direct discretization algorithms. The objective is quadratic and described in ~\eqref{eq:obj-q}. }
	\label{fig:quadratics}
\end{figure}

In this section, we implement our method and test its performance with two different objectives. In the first example, we run the direct discretization (DD) algorithms with different integration order $s$ for solving the strongly convex quadratic problem
\begin{align} \label{eq:obj-q}
	f(x) = x^T Diag(\lambda) x.
\end{align}
$\lambda$ is a vector of singular values. The condition number $\frac{\lambda_{max}}{\lambda_{min}} = 500$. Initial guess $x_0$ is all ones. We compare our methods with GD and NAG algorithms. We use the standard parameter choice for GD and NAG by substituting in the value for smoothness constant $L$ and strongly convex constant $\mu$. For our discretization algorithm, we choose the step size of the format $h = 10^{z},$ where z is the largest integer such that the discretization is stable. The convergence trajectory is shown in Figure~\ref{fig:quadratics}.

\begin{figure}[t]
	\centering
	\vspace{-0.5cm}
	\includegraphics[width=0.38\textwidth]{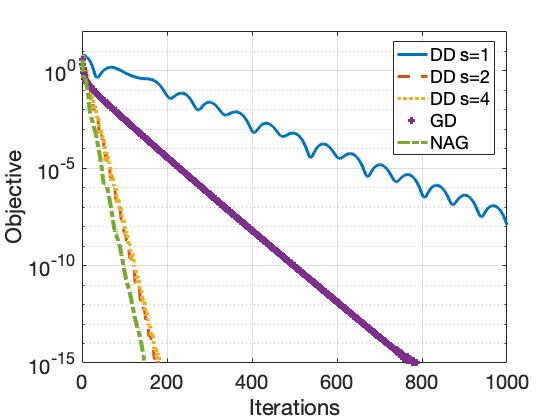}
	\caption{Convergence paths for GD, NAG and direct discretization algorithms. The objective is a logistic regression loss described in ~\eqref{eq:obj-l}. }
	\label{fig:logistics}
	\vspace{-0.5cm}
\end{figure}

In the second example, we tested our example for regularized logistics regression,
\begin{align} \label{eq:obj-l}
	f(w) = \sum_{i=1}^N  \log(1 + e^{-y_i x_i^T w}) + \frac{\gamma}{2} \|w\|^2.
\end{align}
We generate data with a two-cluster mixture of Gaussian model. Each cluster belongs to one class. Particularly, we set the margin large enough such that the data is linearly separable. Then the function is at least $\gamma-$strongly convex. We set the parameter for NAG using $\mu = \gamma$. We then scan the value for smoothness constant $L= 10^{z},$ where z is the largest integer such that the discretization is stable for both GD and NAG. The results are shown in Figure~\ref{fig:logistics}.

Both experiments show similar results and confirm that direct discretization can be faster than gradient descent when the integration order $s > 1$. We also notice that for $s=4$, direct discretization can be as fast as NAG or sometimes even faster. We provide a conjecture for this. First, we use the sub-optimal NAG parameter choice where step size is $\frac{1}{L}$. Second, in the analysis, our convergence rate critically relies on the trajectory of the discretization follows closely to the true solution. However, this requirement is sufficient but not necessary. In fact, by looking at the figure, we notice that the suboptimality of the sequence generated by high-order discretization techinique does not oscillate, which suggests that they are not following the solution to the ODE. Hence, our convergence guarantee is conservative.

	%%% TeX-master: "../root"
\section{Conclusion}
In this work, we showed that by directly discretizating the second order heavy ball ODE, we can obtain optimization algorithms that are provably faster than gradient descent. However, the acceleration can only be shown to exist in a local region surrounding the optimal solution.  This locality is different from traditional locality assumption, because as the function nonlinearity $L$ increase, the neighborhood grows larger. We believe this phenomenon  results from the fact that our proof does not rely on the convexity of the function.  Whether convexity can lead to faster rate in direct discretization remains an interesting research question.

\section{Acknowledgement}
The authors thank Aryan Mokhtari, Cesar Uribe and Juncal Arbelaiz Mugica for helpful discussions.

	\bibliographystyle{IEEEtran}
	\bibliography{root}
\end{document}